\documentclass[12pt]{amsart}
\usepackage{epsfig}

\topskip  \newtheorem{theorem}{Theorem}

\newtheorem{lemma}[theorem]{Lemma}
\newtheorem{corollary}[theorem]{Corollary}

\newcommand\beq{\begin{equation}}
\newcommand\eeq{\end{equation}}
\newcommand\bea{\begin{eqnarray}}
\newcommand\eea{\end{eqnarray}}
\newcommand\bce{\begin{center}}
\newcommand\ece{\end{center}}
\newcommand\ben{\begin{enumerate}}
\newcommand\een{\end{enumerate}}
\newcommand\bit{\begin{itemize}}
\newcommand\eit{\end{itemize}}
\newcommand\nn{\nonumber}

\newcommand\bt{\begin{tabular}}
\newcommand\et{\end{tabular}}
\renewcommand\S{S}
\newcommand\D{{\mathcal D}}
\newcommand\E{{\mathcal E}}
\newcommand\M{{\mathcal M}}
\newcommand\fp{\mathrm{fp}}

\newcommand\wt{\widetilde}
\newcommand\gt{\mathrm{gt}}
\newcommand\hor{\mathrm{hor}}
\newcommand\kra{\varphi}
\newcommand\bmot{\Theta}
\newcommand\lis{\mathrm{lis}}
\newcommand\lds{\mathrm{lds}}
\def\mn{\mbox{-}}
\def\pa{132}
\def\pb{1\mn23}
\def\SS{S}
\def\ttx{{\left(\frac{1-x}{2x}\right)}}

\def\MP{\frak M}
\def\G{N}
\def\vr{\emptyset}

\begin{document}
\title[Restricted Motzkin permutations]{Restricted Motzkin permutations, Motzkin paths, continued fractions, and Chebyshev polynomials}
\maketitle

\begin{center}
{Sergi Elizalde \\
Department of Mathematics \\
MIT, Cambridge, MA 02139\\[4pt]
sergi@math.mit.edu\\[14pt]

Toufik Mansour\\
Department of Mathematics\\
Haifa University\\
31905 Haifa, Israel\\[4pt]
toufik@math.haifa.ac.il}
\end{center}

%------------------------------------------------------------------
\section*{Abstract}
We say that a permutation $\pi$ is a {\em Motzkin permutation} if
it avoids $132$ and there do not exist $a<b$ such that
$\pi_a<\pi_b<\pi_{b+1}$. We study the distribution of several
statistics in Motzkin permutations, including the length of the
longest increasing and decreasing subsequences and the number of
rises and descents. We also enumerate Motzkin permutations with
additional restrictions, and study the distribution of occurrences
of fairly general patterns in this class of permutations.

\noindent{2000 Mathematics Subject Classification}: Primary 05A05,
05A15; Secondary 30B70, 42C05
%------------------------------------------------------------------
\section{Introduction}

\subsection{Background}\label{background}
Let $\alpha\in\SS_n$ and $\tau\in \SS_k$ be two permutations. We
say that $\alpha$ {\it contains\/} $\tau$ if there exists a
subsequence $1\leq i_1<i_2<\dots<i_k\leq n$ such that
$(\alpha_{i_1}, \dots,\alpha_{i_k})$ is order-isomorphic to
$\tau$; in such a context $\tau$ is usually called a {\it
pattern\/}. We say that $\alpha$ {\it avoids\/} $\tau$, or is
$\tau$-{\it avoiding\/}, if such a subsequence does not exist. The
set of all $\tau$-avoiding permutations in $\SS_n$ is denoted
$\SS_n(\tau)$. For an arbitrary finite collection of patterns $T$,
we say that $\alpha$ avoids $T$ if $\alpha$ avoids any $\tau\in
T$; the corresponding subset of $\SS_n$ is denoted $\SS_n(T)$.

While the case of permutations avoiding a single pattern has
attracted much attention, the case of multiple pattern avoidance
remains less investigated. In particular, it is natural, as the next
step, to consider permutations avoiding pairs of patterns $\tau_1$,
$\tau_2$. This problem was solved completely for $\tau_1,\tau_2\in
\SS_3$ (see \cite{SS}), for $\tau_1\in \SS_3$ and $\tau_2\in \SS_4$
(see \cite{W}), and for $\tau_1,\tau_2\in \SS_4$ (see \cite{B1,Km}
and references therein). Several recent papers
\cite{CW,MV1,Kra01,MV2,MV3,MV4} deal with the case $\tau_1\in
\SS_3$, $\tau_2\in \SS_k$ for various pairs $\tau_1,\tau_2$. Another
natural question is to study permutations avoiding $\tau_1$ and
containing $\tau_2$ exactly $t$ times. Such a problem for certain
$\tau_1,\tau_2\in \SS_3$ and $t=1$ was investigated in \cite{R}, and
for certain $\tau_1\in \SS_3$, $\tau_2\in \SS_k$ in
\cite{RWZ,MV1,Kra01}. Most results in these papers are expressed in
terms of Catalan numbers, Chebyshev polynomials, and continued
fractions.

In \cite{BS} Babson and Steingr\'{\i}msson introduced {\em
generalized patterns} that allow the requirement that two adjacent
letters in a pattern must be adjacent in the permutation. In this
context, we write a classical pattern with dashes between any two
adjacent letters of the pattern (for example, $1423$ as
$1\mn4\mn2\mn3$). If we omit the dash between two letters, we mean
that for it to be an occurrence in a permutation $\pi$, the
corresponding letters of $\pi$ have to be adjacent. For example, in
an occurrence of the pattern $12\mn3\mn4$ in a permutation $\pi$,
the letters in $\pi$ that correspond to $1$ and $2$ are adjacent.
For instance, the permutation $\pi=3542617$ has only one occurrence
of the pattern $12\mn3\mn4$, namely the subsequence $3567$, whereas
$\pi$ has two occurrences of the pattern $1\mn2\mn3\mn4$, namely the
subsequences $3567$ and $3467$. Claesson \cite{C} completed the
enumeration of permutations avoiding any single $3$-letter
generalized pattern with exactly one adjacent pair of letters.
Elizalde and Noy \cite{EliNoy} studied some cases of avoidance of
patterns where all letters have to occur in consecutive positions.
Claesson and Mansour \cite{CM} (see also \cite{Mg1,Mg2,Mg3})
presented a complete solution for the number of permutations
avoiding any pair of $3$-letter generalized patterns with exactly
one adjacent pair of letters. Besides, Kitaev \cite{Ki} investigated
simultaneous avoidance of two or more $3$-letter generalized
patterns without internal dashes.

A remark about notation: throughout the paper, a pattern represented with no
dashes will always denote a classical pattern (i.e., with no requirement about
elements being consecutive). All the generalized patterns that we will consider
will have at least one dash.

\subsection{Preliminaries}
{\it Catalan numbers} are defined by
$C_n=\frac{1}{n+1}\binom{2n}{n}$ for all $n\geq0$. The generating
function for the Catalan numbers is given by
$C(x)=\frac{1-\sqrt{1-4x}}{2x}$.

{\it Chebyshev polynomials of the second kind\/} (in what follows
just Chebyshev polynomials) are defined by
$U_r(\cos\theta)=\frac{\sin(r+1)\theta}{\sin\theta}$ for $r\geq0$.
Clearly, $U_r(t)$ is a polynomial of degree $r$ in $t$ with integer
coefficients, which satisfies the following recurrence:
\begin{equation}
U_0(t)=1,\ U_1(t)=2t,\ \mbox{and}\ U_r(t)=2tU_{r-1}(t)-U_{r-2}(t)\
\mbox{for all}\ r\geq2.\label{reccheb}
\end{equation}
The same recurrence is used to define $U_r(t)$ for $r<0$ (for
example, $U_{-1}(t)=0$ and $U_{-2}(t)=-1$). Chebyshev polynomials
were invented for the needs of approximation theory, but are also
widely used in various other branches of mathematics, including
algebra, combinatorics, and number theory (see \cite{Ri}). The
relation between restricted permutations and Chebyshev polynomials
was discovered by Chow and West in \cite{CW}, and later was further
studied by Mansour and Vainshtein \cite{MV1,MV2,MV3,MV4}, and
Krattenthaler \cite{Kra01}. %%
%%These results related to a rational function
%%\begin{equation}
%%R_k(x)=\frac{U_{k-1}\tx}{\sqrt{x}U_k\tx}, \label{rcheb}
%%\end{equation}
%%for all $k\geq 1$. It is easy to see that for any $k$, $R_k(x)$ is
%%rational in $x$ and satisfies
%%\begin{equation}R_k(x)=\frac{1}{1-xR_{k-1}(x)},\label{rk}\end{equation}
%%for all $k\geq1$.

Recall that a {\it Dyck path} of length $2n$ is a lattice path in
$\mathbb{Z}^2$ between $(0,0)$ and $(2n,0)$ consisting of up-steps
$(1,1)$ and down-steps $(1,-1)$ which never goes below the
$x$-axis. Denote by $\D_n$ the set of Dyck paths of length $2n$,
and by $\D=\bigcup_{n\geq0}\D_n$ the class of all Dyck paths. If
$D\in\D_n$, we will write $|D|=n$. Recall that a \emph{Motzkin
path} of length $n$ is a lattice path in $\mathbb{Z}^2$ between
$(0,0)$ and $(n,0)$ consisting of up-steps $(1,1)$, down-steps
$(1,-1)$ and horizontal steps $(1,0)$ which never goes below the
$x$-axis. Denote by $\M_n$ the set of {\it Motzkin paths} with $n$
steps, and let $\M=\bigcup_{n\geq0}\M_n$. We will write $|M|=n$ if
$M\in\M_n$. Sometimes it will be convenient to encode each up-step
by a letter $u$, each down-step by $d$, and each horizontal step
by $h$. Denote by $M_n=|\M_n|$ the $n$-th {\it Motzkin number}.
The generating function for these numbers is $M(x)=\frac{1 - x -
\sqrt{1-2x-3x^2}}{2x^2}$.

Define a {\em Motzkin permutation} $\pi$ to be a $132$-avoiding
permutation in which there do not exist indices $a<b$ such that
$\pi_a<\pi_b<\pi_{b+1}$. Otherwise, if such indices exist,
$\pi_a,\pi_b,\pi_{b+1}$ is called an occurrence of the pattern $\pb$
(for instance, see \cite{C}). For example, there are exactly $4$
Motzkin permutations of length $3$, namely, $213$, $231$, $312$, and
$321$. We denote the set of all Motzkin permutations in $S_n$ by
$\MP_n$. The main reason for the term ``Motzkin permutation'' is
that $|\MP_n|=M_n$, as we will see in Section~\ref{sec:bmot}.

It follows from the definition that the set $\MP_n$ is the same as
the set of $132$-avoiding permutations $\pi\in S_n$ where there is
no $a$ such that $\pi_a<\pi_{a+1}<\pi_{a+2}$. Indeed, assume that
$\pi\in\S_n(132)$ has an occurrence of $\pb$, say
$\pi_a<\pi_b<\pi_{b+1}$ with $a<b$. Now, if $\pi_{b-1}>\pi_b$, then
$\pi$ would have an occurrence of $132$, namely
$\pi_a\pi_{b-1}\pi_{b+1}$. Therefore, $\pi_{b-1}<\pi_b<\pi_{b+1}$,
so $\pi$ has three consecutive increasing elements.

For any subset $A\in S_n$ and any pattern $\alpha$, define
$A(\alpha):=A\cap S_n(\alpha)$. For example, $\MP_n(\alpha)$
denotes the set of Motzkin permutations of length $n$ that avoid
$\alpha$.

%%An {\em involution} $\pi$ is a permutation such that
%%$\pi=\pi^{-1}$;  let $I_n$ denote the set of all $n$-involutions.
%%Pattern-avoiding involutions have also been linked with other
%%combinatorial objects (see Guibert \cite{GuibertThese} and
%%refrences therein). In particularly, Guibert \cite{GuibertThese}
%%has established bijection between 1-2 trees having $n$ edges and
%%the set $I_n(3412)$ where $I_n(3412)$ is enumerated by the $n$-th
%%Motzkin number. Recently, Egge~\cite{Egge} gave enumerations and
%%generating functions for involutions in $I_n(3412)$ which avoid
%%various sets of additional patterns.

\subsection{Organization of the paper}
In Section~\ref{sec:bmot} we exhibit a bijection between the set
of Motzkin permutations and the set of Motzkin paths. Then we use
it to obtain generating functions of Motzkin permutations with
respect to the length of the longest decreasing and increasing
subsequences together with the number of rises. The section ends
with another application of the bijection, to the enumeration of
fixed points in permutations avoiding simultaneously $231$ and
$32\mn1$.

In Section~\ref{sec:restricted} we consider additional restrictions
on Motzkin permutations. Using a block decomposition, we enumerate
Motzkin permutations avoiding the pattern $12\ldots k$, and we find
the distribution of occurrences of this pattern in Motzkin
permutations. Then we obtain generating functions for Motzkin
permutations avoiding patterns of more general shape. We conclude
the section by considering two classes of generalized patterns (as
described above), and we study its distribution in Motzkin
permutations.

%----------------------------------------------------------------
\section{Bijection $\bmot:\MP_n\longrightarrow\M_n$}\label{sec:bmot}
In this section we establish a bijection $\bmot$ between Motzkin
permutations and Motzkin paths. This bijection allows us to describe
the distribution of some interesting statistics on the set of
Motzkin permutations.

\subsection{The bijection $\bmot$}
\begin{figure}[hbt]
\epsfig{file=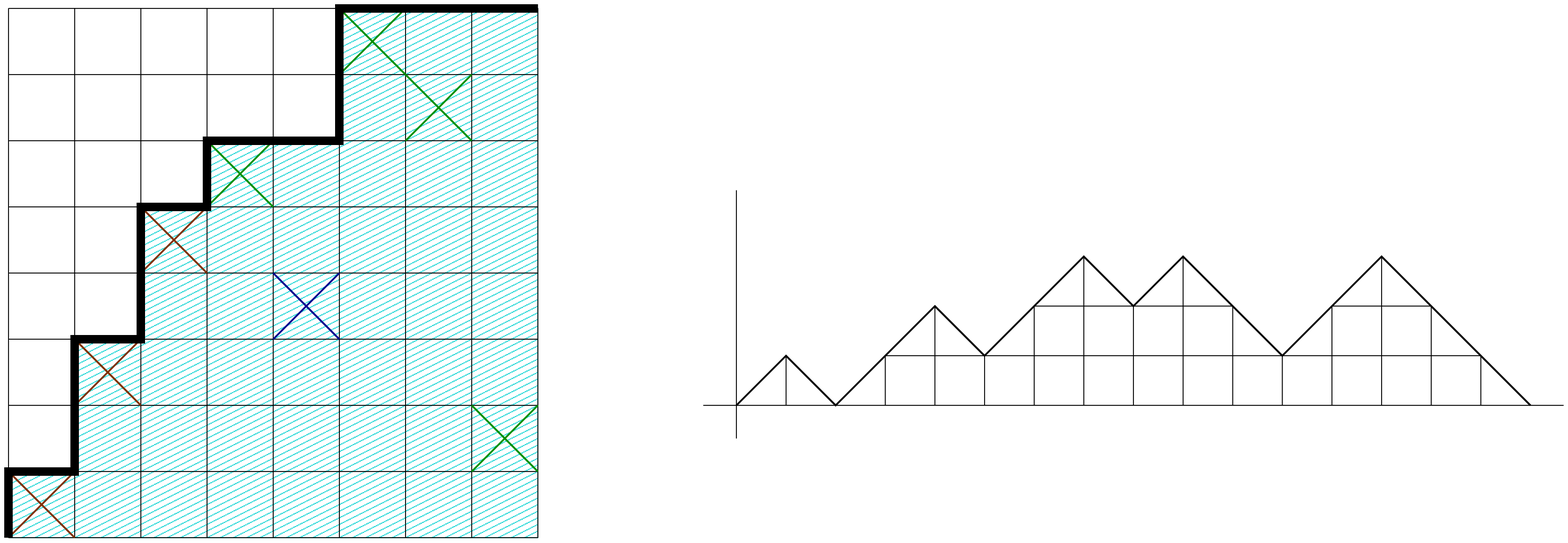,height=3.7cm} \caption{\label{fig:bij2}
The bijection $\kra$.}
\end{figure}

We can give a bijection $\bmot$ between $\MP_n$ and $\M_n$. In order
to do so we use first the following bijection $\kra$ from
$\S_n(132)$ to $\D_n$, which is essentially due to Krattenthaler
\cite{Kra01}, and also described independently by Fulmek
\cite{Ful01} and Reifegerste \cite{Rei02}. Consider
$\pi\in\S_n(132)$ given as an $n\times n$ array with crosses in the
squares $(i,\pi_i)$. Take the path with \emph{up} and \emph{right}
steps that goes from the lower-left corner to the upper-right
corner, leaving all the crosses to the right, and staying always as
close to the diagonal connecting these two corners as possible. Then
$\kra(\pi)$ is the Dyck path obtained from this path by reading an
up-step every time the path goes up and a down-step every time it
goes right. Figure~\ref{fig:bij2} shows an example when
$\pi=67435281$.

There is an easy way to recover $\pi$ from $\kra(\pi)$. Assume we
are given the path from the lower-left corner to the upper-right
corner of the array. Row by row, put a cross in the leftmost square
to the right of this path such that there is exactly one cross in
each column. This gives us $\pi$ back.

One can see that $\pi\in\S_n(132)$ avoids $\pb$ if and only if the
Dyck path $\kra(\pi)$ does not contain three consecutive up-steps
(a \emph{triple rise}). Indeed, assume that $\kra(\pi)$ has three
consecutive up-steps. Then, the path from the lower-left corner to
the upper-right corner of the array has three consecutive vertical
steps. The crosses in the corresponding three rows give three
consecutive increasing elements in $\pi$ (this follows from the
definition of the inverse of $\kra$), and hence an occurrence of
$\pb$.% (see remark below).

Reciprocally, assume now that $\pi$ has an occurrence of $\pb$.
The path from the lower-left to the upper-right corner of the
array of $\pi$ must have two consecutive vertical steps in the
rows of the crosses corresponding to `2' and `3'. But if
$\kra(\pi)$ has no triple rise, the next step of this path must be
horizontal, and the cross corresponding to `2' must be right below
it. But then all the crosses above this cross are to the right of
it, which contradicts the fact that this was an occurrence of
$\pb$.

%\begin{remark}
%Actually we have shown that if $\kra(\pi)$ contains a triple rise,
%then there exists an index $a$ such that
%$\pi_a<\pi_{a+1}<\pi_{a+2}$. This implies that the set $\MP_n$ is
%the same as the set of $132$-avoiding permutations $\pi\in S_n$
%where there is no $a$ such that $\pi_a<\pi_{a+1}<\pi_{a+2}$.
%\end{remark}

Denote by $\E_n$ the set of Dyck paths of length $2n$ with no
triple rise. We have given a bijection between $\MP_n$ and $\E_n$.
The second step is to exhibit a bijection between $\E_n$ and
$\M_n$, so that $\bmot$ will be defined as the composition of the
two bijections. Given $D\in\E_n$, divide it in $n$ blocks,
splitting after each down-step. Since $D$ has no triple rises,
each block is of one of these three forms: $uud$, $ud$, $d$. From
left to right, transform the blocks according to the rule \bea\nn
uud & \rightarrow & u,\\ \label{rule} ud & \rightarrow & h,\\ \nn
d & \rightarrow & d.\eea We obtain a Motzkin path of length $n$.
This step is clearly a bijection.

Up to reflection of the Motzkin path over a vertical line, $\bmot$
is essentially the same bijection that was given by Claesson
\cite{C} between $\MP_n$ and $\M_n$, using a recursive definition.

\subsection{Statistics in $\MP_n$}
Here we show applications of the bijection $\bmot$ to give
generating functions for several statistics on Motzkin permutations.
For a permutation $\pi$, denote by $\lis(\pi)$ and $\lds(\pi)$
respectively the length of the longest increasing subsequence and
the length of the longest decreasing subsequence of $\pi$. The
following lemma follows from the definitions of the bijections and
from the properties of $\kra$ (see \cite{Kra01}).

\begin{lemma}\label{lemma:bmot}
Let $\pi\in\MP_n$, let $D=\kra(\pi)\in\D_n$, and let
$M=\bmot(\pi)\in\M_n$. We have \ben \item $\lds(\pi) = \#\{$peaks
of $D\} = \#\{$steps $u$ in $M\} + \#\{$steps $h$ in $M\}$,
\item $\lis(\pi) =$ height of $D$ = height of $M + 1$,
\item $\#\{$rises of $\pi\} = \#\{$double rises of $D\} = \#\{$steps $u$ in
$M\}$. \een
\end{lemma}

\begin{theorem}\label{th:lds} The generating function for Motzkin permutations
with respect to the length of the longest decreasing subsequence
and to the number of rises is \bea\nn
A(v,y,x):=\sum_{n\ge0}\sum_{\pi\in\MP_n}v^{\lds(\pi)}y^{\#\{\mathrm{rises\
of\ }\pi\}}x^n= \frac{1-vx-\sqrt{1-2vx+(v^2-4vy)x^2}}{2v y
x^2}.\eea Moreover,
$$A(v,y,x)=\sum_{n\geq0}\sum_{m\geq0}\frac{1}{n+1}\binom{2n}{n}\binom{m+2n}{2n}x^{m+2n}v^{m+n}y^n.$$
\end{theorem}
\begin{proof}
By Lemma~\ref{lemma:bmot}, we can express $A$ as
$$A(v,y,x)=\sum_{M\in\M}v^{\#\{\mathrm{steps\ } u \mathrm{\ in\ } M\} + \#\{\mathrm{steps\ } h \mathrm{\ in\ } M\}}
y^{\#\{\mathrm{steps\ } u \mathrm{\ in\ } M\}}x^{|M|}.$$ Using the
standard decomposition of Motzkin paths, we obtain the following
equation for the generating function $A$. \bea\label{eqA}
A(v,y,x)=1 + v x A(v,y,x) + v y x^2 A^2(v,y,x).\eea Indeed, any
nonempty $M\in\M$ can be written uniquely in one of the following
two forms: \ben \item $M=hM_1$, \item $M=uM_1dM_2$,\een where
$M_1,M_2,M_3$ are arbitrary Motzkin paths. In the first case, the
number of horizontal steps of $hM_1$ is one more than in $M_1$,
the number of up steps is the same, and $|hM_1|=|M_1|+1$, so we
get the term $v x A(v,y,x)$. Similarly, the second case gives the
term $v y x^2 A^2(v,y,x)$. Solving equation~(\ref{eqA}) we get
that
$$A(v,y,x)= \frac{1-vx-\sqrt{1-2vx+(v^2-4vy)x^2}}{2v y x^2}=\frac{1}{1-vx}C\left(\frac{vyx^2}{(1-vx)^2}\right),$$
where $C(t)=\frac{1-\sqrt{1-4t}}{2t}$ the generating function for
the Catalan numbers. Thus,
$$A(v,y,x)= \sum_{n\geq0}\frac{1}{n+1}\binom{2n}{n}\frac{y^nx^{2n}v^n}{(1-vx)^{2n+1}}=
\sum_{n\geq0}\sum_{m\geq0}\frac{1}{n+1}\binom{2n}{n}\binom{m+2n}{2n}x^{m+2n}v^{m+n}y^n.$$
\end{proof}

\begin{theorem}\label{th:lislds} For $k>0$, let
$$B_k(v,y,x):=\sum_{n\ge0}\sum_{\pi\in\MP_n(12\ldots(k+1))}v^{\lds(\pi)}
y^{\#\{\mathrm{rises\ of\ }\pi\}}x^n$$ be the generating function
for Motzkin permutations avoiding $12\ldots(k+1)$ with respect to
the length of the longest decreasing subsequence and to the number
of rises. Then we have the recurrence
$$B_k(v,y,x)=\frac{1}{1-vx-vyx^2B_{k-1}(v,y,x)},$$ with
$B_1(v,y,x)=\frac{1}{1-vx}$. Thus, $B_k$ can be expressed as
$$B_k(v,y,x)=\frac{1}{1-vx-\dfrac{vyx^2}{1-vx-\dfrac{vyx^2}{\dfrac{\ddots}{1-vx-\dfrac{vyx^2}{1-vx}}}}},$$
where the fraction has $k$ levels, or in terms of Chebyshev
polynomials of the second kind, as
$$B_k(v,y,x)=\frac{U_{k-1}\left(\frac{1-vx}{2x\sqrt{vy}}\right)}{x\sqrt{vy}U_k\left(\frac{1-vx}{2x\sqrt{vy}}\right)}.$$
\end{theorem}
\begin{proof}
The condition that $\pi$ avoids $12\ldots(k+1)$ is equivalent to
the condition $\lis(\pi)\le k$. By Lemma~\ref{lemma:bmot},
permutations in $\MP_n$ satisfying this condition are mapped by
$\bmot$ to Motzkin paths of height strictly less than $k$. Thus,
we can express $B_k$ as
$$B_k(v,y,x)=\underset{\mathrm{of\ height}<k}{\sum_{M\in\M}}v^{\#\{ \mathrm{steps\ } u \mathrm{\ in\ } M\} + \#\{\mathrm{steps\ } h \mathrm{\ in\ }M\}}
y^{\#\{\mathrm{steps\ } u \mathrm{\ in\ } M\}}x^{|M|}.$$ The
continued fraction follows now from \cite{Fla80}. Alternatively, we
can use again the standard decomposition of Motzkin paths, for
$k>1$. In the first of the above cases, the height of $hM_1$ is the
same as the height of $M_1$. However, in the second case, in order
for the height of $uM_2dM_3$ to be less than $k$, the height of
$M_2$ has to be less than $k-1$. So we obtain the equation
$$B_k(v,y,x)=1 + v x B_k(v,y,x) + v y x^2 B_{k-1}(v,y,x)B_k(v,y,x).$$ For $k=1$,
the path can have only horizontal steps, so we get
$B_1(v,y,x)=\frac{1}{1-vx}$. Now, using the above recurrence and
Equation~\ref{reccheb} we get the desired result.
\end{proof}

\subsection{Fixed points in the reversal of Motzkin permutations}
Here we show another application of $\bmot$. A slight modification
of it will allow us to enumerate fixed points in another class of
pattern-avoiding permutations closely related to Motzkin
permutations. For any $\pi=\pi_1\pi_2\ldots\pi_n\in\S_n$, denote
its reversal by $\pi^R=\pi_n\ldots\pi_2\pi_1$. Let
$\MP_n^R:=\{\pi\in S_n:\pi^R\in\MP_n\}$. In terms of pattern
avoidance, $\MP_n^R$ is the set of permutations that avoid $231$
and $32\mn1$ simultaneously, that is, the set of 231-avoiding
permutations $\pi\in S_n$ where there do not exist $a<b$ such that
$\pi_{a-1}>\pi_a>\pi_b$. Recall that $i$ is called a {\em fixed
point} of $\pi$ if $\pi_i=i$.

\begin{theorem}
The generating function
$\sum_{n\ge0}\sum_{\pi\in\MP_n^R}w^{\mathrm{fp}(\pi)}x^n$ for
permutations avoiding simultaneously $231$ and $32\mn1$ with
respect to to the number of fixed points is
\bea\nn\hspace{12cm}\\
\hspace{0.4cm}\dfrac{1}{1-wx-\dfrac{x^2}{1-x-M_0(w-1)x^2-\dfrac{x^2}{1-x-M_1(w-1)x^3-\frac{x^2}{1-x-M_2(w-1)x^4-\frac{x^2}
{\ddots}}}}},\label{fpMPR}\eea where after the second level, the
coefficient of $(w-1)x^{n+2}$ is the Motzkin number $M_n$.
\end{theorem}
\begin{proof}
We have the following composition of bijections: \bce\bt{ccccccc}
$\MP_n^R$ & $\longleftrightarrow$ & $\MP_n$ &
$\longleftrightarrow$ & $\E_n$ & $\longleftrightarrow$ & $\M_n$ \\
$\pi$ & $\mapsto$ & $\pi^R$ & $\mapsto$ & $\kra(\pi^R)$ &
$\mapsto$ & $\bmot(\pi^R)$ \et\ece The idea of the proof is to
look at how the fixed points of $\pi$ are transformed by each of
these bijections.

We use the definition of tunnel of a Dyck path given
in~\cite{Eli02}, and generalize it to Motzkin paths. A \emph{tunnel}
of $M\in\M$ (resp. $D\in\D$) is a horizontal segment between two
lattice points of the path that intersects $M$ (resp. $D$) only in
these two points, and stays always below the path. Tunnels are in
obvious one-to-one correspondence with decompositions of the path as
$M=XuYdZ$ (resp. $D=XuYdZ$), where $Y\in\M$ (resp. $Y\in\D$). In the
decomposition, the tunnel is the segment that goes from the
beginning of the $u$ to the end of the $d$. Clearly such a
decomposition can be given for each up-step $u$, so the number of
tunnels of a path equals its number of up-steps. The {\it length} of
a tunnel is just its length as a segment, and the {\it height} is
the $y$-coordinate of the segment.

Fixed points of $\pi$ are mapped by the reversal operation to
elements $j$ such that $\pi^R_j=n+1-j$, which in the array of
$\pi^R$ correspond to crosses on the diagonal between the
bottom-left and top-right corners. Each cross in this array
naturally corresponds to a tunnel of the Dyck path $\kra(\pi^R)$,
namely the one determined by the vertical step in the same row as
the cross and the horizontal step in the same column as the cross.
It is not hard to see (and is also shown in~\cite{Eli03}) that
crosses on the diagonal between the bottom-left and top-right
corners correspond in the Dyck path to tunnels $T$ satisfying the
condition $\mathrm{height}(T)+1=\frac{1}{2}\mathrm{length}(T)$.

The next step is to see how these tunnels are transformed by the
bijection from $\E_n$ to $\M_n$. Tunnels of height 0 and length 2
in the Dyck path $D:=\kra(\pi^R)$ are just hills $ud$ landing on
the $x$-axis. By the rule (\ref{rule}) they are mapped to
horizontal steps at height 0 in the Motzkin path
$M:=\bmot(\pi^R)$. Assume now that $k\ge 1$. A tunnel $T$ of
height $k$ and length $2(k+1)$ in $D$ corresponds to a
decomposition $D=XuYdZ$ where $X$ ends at height $k$ and
$Y\in\D_{2k}$. Note that $Y$ has to begin with an up-step (since
it is a nonempty Dyck path) followed by a down-step, otherwise $D$
would have a triple rise. Thus, we can write $D=XuudY'dZ$ where
$Y'\in\D_{2(k-1)}$. When we apply to $D$ the bijection given by
rule (\ref{rule}), $X$ is mapped to an initial segment $\wt{X}$ of
a Motzkin path ending at height $k$, $uud$ is mapped to $u$, $Y'$
is mapped to a Motzkin path $\wt{Y'}\in\M_{k-1}$ of length $k-1$,
the $d$ following $Y'$ is mapped to $d$ (since it is preceded by
another $d$), and $Z$ is mapped to a final segment $\wt{Z}$ of a
Motzkin path going from height $k$ to the $x$-axis. Thus, we have
that $M=\wt{X}u\wt{Y'}d\wt{Z}$. It follows that tunnels $T$ of $D$
satisfying $\mathrm{height}(T)+1=\frac{1}{2}\mathrm{length}(T)$
are transformed by the bijection into tunnels $\wt{T}$ of $M$
satisfying $\mathrm{height}(\wt{T})+1=\mathrm{length}(\wt{T})$. We
will call {\it good} tunnels the tunnels of $M$ satisfying this
last condition. It remains to show that the generating function
for Motzkin paths where $w$ marks the number of good tunnels plus
the number of horizontal steps at height 0, and $x$ marks the
length of the path, is given by (\ref{fpMPR}).

To do this we imitate the technique used in \cite{Eli03} to
enumerate fixed points in 231-avoiding permutations. We will
separate good tunnels according to their height. It is important to
notice that if a good tunnel of $M$ corresponds to a decomposition
$M=XuYdZ$, then $M$ has no good tunnels inside the part given by
$Y$. In other words, the orthogonal projections on the $x$-axis of
all the good tunnels of a given Motzkin path are disjoint. Clearly,
they are also disjoint from horizontal steps at height 0. Using this
observation, one can apply directly the results in \cite{Fla80} to
give a continued fraction expression for our generating function.
However, for the sake of completeness we will explain here how to
obtain this expression.

For every $k\ge1$, let $\gt_k(M)$ be the number of tunnels of $M$
of height $k$ and length $k+1$. Let $\hor(M)$ be the number of
horizontal steps at height 0. We have seen that for
$\pi\in\MP_n^R$,
$\fp(\pi)=\hor(\bmot(\pi^R))+\sum_{k\ge1}\gt_k(\bmot(\pi^R))$. We
will show now that for every $k\ge1$, the generating function for
Motzkin paths where $w$ marks the statistic
$\hor(M)+\gt_1(M)+\cdots+\gt_{k-1}(M)$ is given by the continued
fraction~(\ref{fpMPR}) truncated at level $k$, with the $(k+1)$-st
level replaced with $M(x)$.

A Motzkin path $M$ can be written uniquely as a sequence of
horizontal steps $h$ and elevated Motzkin paths $uM'd$, where
$M'\in\M$. In terms of the generating function
$M(x)=\sum_{M\in\M}x^{|M|}$, this translates into the equation
$M(x)=\frac{1}{1-x-x^2M(x)}$. The generating function where $w$
marks horizontal steps at height 0 is just
$$\sum_{M\in\M}w^{\hor(M)}x^{|M|}=\frac{1}{1-wx-x^2M(x)}.$$

If we want $w$ to mark also good tunnels at height 1, each $M'$
from the elevated paths above has to be decomposed as a sequence
of horizontal steps and elevated Motzkin paths $uM''d$. In this
decomposition, a tunnel of height 1 and length 2 is produced by
each empty $M''$, so we have
\beq\sum_{M\in\M}w^{\hor(M)+\gt_1(M)}x^{|M|}=\frac{1}{1-wx-\dfrac{x^2}{1-x-x^2[w-1+M(x)]}}.\label{height1}\eeq
Indeed, the $M_0(=1)$ possible empty paths $M''$ have to be
accounted as $w$, not as 1.

Let us now enumerate simultaneously horizontal steps at height 0
and good tunnels at heights 1 and 2. We can rewrite
(\ref{height1}) as
$$\frac{1}{1-wx-\dfrac{x^2}{1-x-x^2\left[w-1+\dfrac{1}{1-x-x^2M(x)}\right]}}.$$
Combinatorially, this corresponds to expressing each $M''$ as a
sequence of horizontal steps and elevated paths $uM'''d$, where
$M'''\in\M$. Notice that since $uM'''d$ starts at height 2, a
tunnel of height 2 and length 3 is created whenever $M'''\in\M_1$.
Thus, if we want $w$ to mark also these tunnels, such an $M'''$
has to be accounted as $wx$, not $x$. The corresponding generating
function is
$$\begin{array}{l}
\sum_{M\in\M}w^{\hor(M)+\gt_1(M)+\gt_2(M)}x^{|M|}\\
\qquad\qquad\qquad=\dfrac{1}{1-wx-\dfrac{x^2}{1-x-x^2\left[w-1+
\dfrac{1}{1-x-x^2[(w-1)x+M(x)]}\right]}}. \end{array}$$

Now it is clear how iterating this process indefinitely we obtain
the continued fraction (\ref{fpMPR}). From the generating function
where $w$ marks $\hor(M)+\gt_1(M)+\cdots+\gt_k(M)$, we can obtain
the one where $w$ marks $\hor(M)+\gt_1(M)+\cdots+\gt_{k+1}(M)$ by
replacing the $M(x)$ at the lowest level with
$$\frac{1}{1-x-x^2[M_k(w-1)x^k+M(x)]},$$ to account for tunnels of
height $k$ and length $k+1$, which in the decomposition correspond
to elevated Motzkin paths at height $k$.
\end{proof}

\section{Restricted Motzkin permutations}\label{sec:restricted}
In this section we consider those Motzkin permutations in $\MP_n$
that avoid an arbitrary pattern $\tau$. More generally, we
enumerate Motzkin permutations according to the number of
occurrences of $\tau$. Subsection~\ref{ss:12k} deals with the
increasing pattern $\tau=12\ldots k$. In
Subsection~\ref{ss:subpattern} we show that if $\tau$ has a
certain form, we can express the generating function for
$\tau$-avoiding Motzkin permutations in terms of the the
corresponding generating functions for some subpatterns of $\tau$.
Finally, Subsection~\ref{ss:gener} studies the case of the
generalized patterns $12\mn3\mn\ldots\mn k$ and
$21\mn3\mn\ldots\mn k$.

We begin by introducing some notation. Let $M_\tau(n)$ be the number
of Motzkin permutations in $\MP_n(\tau)$, and let
$\G_\tau(x)=\sum_{n\geq0}M_\tau(n)x^n$ be the corresponding
generating function.

\begin{center}
\begin{figure}[h]
\epsfxsize=2.4in \epsffile{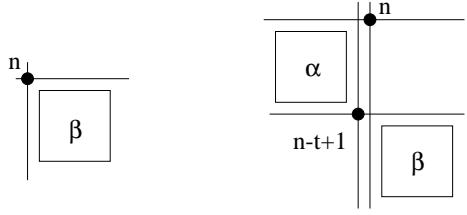} \caption{The block
decomposition for $\pi\in\MP_n$.} \label{avmp}
\end{figure}
\end{center}

Let $\pi\in\MP_n$. Using the block decomposition approach
(see~\cite{MV4}), we have two possible block decompositions of
$\pi$, as shown in Figure~\ref{avmp}. These decompositions are
described in Lemma~\ref{pmp}, which is the basis for all the results
in this section.

\begin{lemma}\label{pmp} Let
$\pi\in\MP_n$. Then one of the following holds:

{\rm(i)} $\pi=(n,\beta)$ where $\beta\in\MP_{n-1}$,

{\rm(ii)} there exists $t$, $2\leq t\leq n$, such that
$\pi=(\alpha,n-t+1,n,\beta)$, where
$$(\alpha_1-(n-t+1),\ldots,\alpha_{t-2}-(n-t+1))\in\MP_{t-2}\mbox{
and }\beta\in\MP_{n-t}.$$
\end{lemma}
\begin{proof}
Given $\pi\in\MP_n$, take $j$ so that $\pi_j=n$. Then
$\pi=(\pi',n,\pi'')$, and the condition that $\pi$ avoids $\pa$ is
equivalent to $\pi'$ being a permutation of the numbers $n-j+1,
n-j+2,\ldots,n-1$, $\pi''$ being a permutation of the numbers
$1,2,\dots,n-j$, and both $\pi'$ and $\pi''$ being $\pa$-avoiding.
On the other hand, it is easy to see that if $\pi'$ is nonempty,
then $\pi$ avoids $\pb$ if and only if the minimal entry of $\pi'$
is adjacent to $n$, and both $\pi'$ and $\pi''$ avoid $\pb$.
Therefore, $\pi$ avoids $\pa$ and $\pb$ if and only if either (i)
or (ii) hold.
\end{proof}

\subsection{The pattern $\tau=12\ldots k$}\label{ss:12k}
From Theorem~\ref{th:lislds} we get the following expression for
$\G_\tau$:
  $$\G_{12\ldots k}(x)=\frac{U_{k-2}\ttx}{xU_{k-1}\ttx}.$$
This result can also be easily proved using the block
decomposition given in Lemma~\ref{pmp}. Now we turn our attention
to analogues of \cite[Theorem~1]{BCS}. Let $\G(x_1,x_2,\ldots)$ be
the generating function
$$\sum_{n\geq0}\sum_{\pi\in\MP_n}\prod_{j\geq1} x_j^{12\ldots
j(\pi)},$$ where $12\ldots j(\pi)$ is the number of occurrences of
the pattern $12\ldots j$ in $\pi$.

\begin{theorem}\label{thac}
The generating function
$\sum\limits_{n\geq0}\sum\limits_{\pi\in\MP_n} \prod\limits_{j\geq
1} x_j^{12\ldots j}(\pi)$ is given by the following continued
fraction:
$$
\dfrac{1}{1-x_1-\dfrac{x_1^2x_2}{1-x_1x_2-\dfrac{x_1^2x_2^3x_3}{1-x_1x_2^2x_3
-\dfrac{x_1^2x_2^5x_3^4x_4}{\ddots}}}},$$ in which the $n$-th
numerator is
$\prod\limits_{i=1}^{n+1}x_i^{\binom{n}{i-1}+\binom{n-1}{i-1}}$
and the $n$-th denominator is
$\prod\limits_{i=1}^{n}x_i^{\binom{n-1}{i-1}}$.
\end{theorem}
\begin{proof}
By Lemma~\ref{pmp}, we have two possibilities for the block
decomposition of an arbitrary Motzkin permutation $\pi\in\MP_n$.
Let us write an equation for $\G(x_1,x_2,\ldots)$. The
contribution of the first decomposition is
$x_1\G(x_1,x_2,\ldots)$, and the second decomposition gives
$x_1^2x_2\G(x_1x_2,x_2x_3,\ldots)\G(x_1,x_2,\ldots)$. Therefore,
$$\G(x_1,x_2,\ldots)=1+x_1\G(x_1,x_2,\ldots)+x_1^2x_2
\G(x_1x_2,x_2x_3,\ldots)\G(x_1,x_2,\ldots),$$ where $1$ is the
contribution of the empty Motzkin permutation. The theorem follows
now by induction.
\end{proof}

\subsubsection{Counting occurrences of the pattern $12\ldots k$ in a Motzkin
permutation} Using Theorem~\ref{thac} we can enumerate occurrences
of the pattern $12\ldots k$ in Motzkin permutations.

\begin{theorem}\label{thacapp1}
Fix $k\geq2$. The generating function for the number of Motzkin
permutations which contain $12\ldots k$ exactly $r$ times is given
by
   $$\frac{\left(U_{k-2}\ttx-xU_{k-3}\ttx\right)^{r-1}}{U_{k-1}^{r+1}\ttx},$$
for all $r=1,2,\ldots,k$.
\end{theorem}
\begin{proof}
Let $x_1=x$, $x_k=y$, and $x_j=1$ for all $j\neq 1,k$. Let
$G_k(x,y)$ be the function obtained from $\G(x_1,x_2,\ldots)$
after this substitution. Theorem~\ref{thac} gives
$$G_k(x,y)=\dfrac{1}{1-x-\dfrac{x^2}{1-x-\dfrac{x^2}{\ddots-\dfrac{\ddots}{1-x-\dfrac{x^2y}{1-xy-\dfrac{x^2y^{k+1}}{\ddots}}}}}}.$$
So, $G_k(x,y)$ can be expressed as follows. For all $k\geq2$,
$$G_k(x,y)=\frac{1}{1-x-x^2G_{k-1}(x,y)},$$ and there exists a
continued fraction $H(x,y)$ such that
$G_1(x,y)=\frac{y}{1-xy-y^{k+1}H(x,y)}$. Now, using induction on
$k$ together with (\ref{reccheb}) we get that there exists a
formal power series $J(x,y)$ such that
$$G_k(x,y)=\frac{U_{k-2}\ttx-\left(U_{k-3}\ttx-xU_{k-4}\ttx\right)y}
{xU_{k-1}\ttx-x\left(U_{k-2}\ttx-xU_{k-3}\ttx\right)y}+y^{k+1}J(x,y).$$
The series expansion of $G_k(x,y)$ about the point $y=0$ gives
$$\begin{array}{l}
G_k(x,y)=\left[U_{k-2}\ttx-\left(U_{k-3}\ttx-xU_{k-4}\ttx\right)y\right]\\
\qquad\qquad\qquad\qquad\qquad\qquad\qquad\cdot\sum\limits_{r\geq0}\frac{\left(U_{k-2}\ttx-xU_{k-3}\ttx\right)^r}{xU_{k-1}^{r+1}\ttx}y^r+y^{k+1}J(x,y).
\end{array}$$
Hence, by using the identities
$$U_{k}^2(t)-U_{k-1}(t)U_{k+1}(t)=1\mbox{ and }
U_{k}(t)U_{k-1}(t)-U_{k-2}(t)U_{k+1}(t)=2t$$ we get the desired
result.
\end{proof}

\subsubsection{More statistics on Motzkin permutations}
We can use the above theorem to find the generating function for
the number of Motzkin permutations with respect to various
statistics.

For another application of Theorem~\ref{thac}, recall that $i$ is
a {\em free rise} of $\pi$ if there exists $j$ such that
$\pi_i<\pi_j$. We denote the number of free rises of $\pi$ by
$fr(\pi)$. Using Theorem~\ref{thac} for $x_1=x$, $x_2=q$, and
$x_j=1$ for $j\geq3$, we get the following result.

\begin{corollary}\label{thst1}
The generating function
$\sum_{n\geq0}\sum_{\pi\in\MP_n}x^nq^{fr(\pi)}$  is given by the
following continued fraction:
$$
\dfrac{1}{1-x-\dfrac{x^2q}{1-xq-\dfrac{x^2q^3}{1-xq^2
-\dfrac{x^2q^5}{\ddots}}}},$$ in which the $n$-th numerator is
$x^2q^{2n-1}$ and the $n$-th denominator is $xq^{n-1}$.
\end{corollary}

For our next application, recall that $\pi_j$ is a {\em
left-to-right maximum} of a permutation $\pi$ if $\pi_i<\pi_j$ for
all $i<j$. We denote the number of left-to-right maxima of $\pi$
by $lrm(\pi)$.

\begin{corollary}\label{thst2}
The generating function
$\sum_{n\geq0}\sum_{\pi\in\MP_n}x^nq^{lrm(\pi)}$  is given by the
following continued fraction:
$$
\dfrac{1}{1-xq-\dfrac{x^2q}{1-x-\dfrac{x^2}{1-x
-\dfrac{x^2}{\ddots}}}}.$$ Moreover,
$$\sum_{n\geq0}\sum_{\pi\in\MP_n}x^nq^{lrm(\pi)}=\sum_{m\geq0}x^m(1+xM(x))^mq^m.$$
\end{corollary}
\begin{proof}
Using Theorem~\ref{thac} for $x_1=xq$, and
$x_{2j}=x_{2j+1}^{-1}=q^{-1}$ for $j\geq1$, together with
\cite[Proposition~5]{BCS} we get the first equation as claimed. The
second equation follows from the fact that the continued fraction
$$\dfrac{1}{1-x-\dfrac{x^2}{1-x -\dfrac{x^2}{\ddots}}}$$
is given by the generating function for the Motzkin numbers,
namely $M(x)$.
\end{proof}

\subsection{General restriction}\label{ss:subpattern}
 Let us find the generating function for those
Motzkin permutations which avoid $\tau$ in terms of the generating
function for Motzkin permutations avoiding $\rho$, where $\rho$ is a
permutation obtained by removing some entries from $\tau$. The next
theorem is analogous to the result for $123$-avoiding permutations
that appears in \cite[Theorem 9]{Kra01}.

\begin{theorem}\label{thf0}
Let $k\geq4$, $\tau=(\rho',1,k)\in\MP_k$, and let
$\rho\in\MP_{k-2}$ be the permutation obtained by decreasing each
entry of $\rho'$ by 1. Then
$$\G_{\tau}(x)=\frac{1}{1-x-x^2\G_\rho(x)}.$$
\end{theorem}
\begin{proof}
By Lemma~\ref{pmp}, we have two possibilities for the block
decomposition of a nonempty Motzkin permutation in $\MP_n$. Let us
write an equation for $\G_{\tau}(x)$. The contribution of the first
decomposition is $x\G_{\tau}(x)$, and from the second decomposition
we get $x^2\G_{\rho}(x)\G_{\tau}(x)$. Hence,
$$\G_{\tau}(x)=1+x\G_\tau(x)+x^2\G_{\rho}(x)\G_{\tau}(x),$$
where $1$ corresponds to the empty Motzkin permutation. Solving
the above equation we get the desired result.
\end{proof}

As an extension of \cite[Theorem 9]{Kra01}, let us consider the
case $\tau=23\ldots (k-1)1k$. Theorem~\ref{thf0} for
$\tau=23\ldots (k-1)1k$ ($\rho=12\ldots(k-2)$) gives
$$\G_{23\ldots(k-1)1k}(x)=\frac{1}{1-x-x^2\G_{12\ldots(k-2)}(x)}.$$
Hence, by Theorem~\ref{th:lislds} together with (\ref{reccheb}) we
get
$$\G_{23\ldots(k-1)1k}(x)=\frac{U_{k-3}\ttx}{xU_{k-2}\ttx}.$$
\begin{corollary}\label{cthf0a}
For all $k\geq1$,
$$\G_{k(k+1)(k-1)(k+2)(k-2)(k+3)\ldots1(2k)}(x)=
\frac{U_{k-1}\ttx}{xU_k\ttx},$$ and
$$\G_{(k+1)k(k+2)(k-1)(k+3)\ldots1(2k+1)}(x)
=\frac{U_{k}\ttx+U_{k-1}\ttx}{x\left(U_{k+1}\ttx+U_{k}\ttx\right)}.$$
\end{corollary}
\begin{proof}
Theorem~\ref{thf0} for
$\tau=k(k+1)(k-1)(k+2)(k-2)(k+3)\ldots1(2k)$ gives
    $$\G_{\tau}(x)=\frac{1}{1-x-x^2\G_{(k-1)k(k-2)(k+1)(k-3)(k+2)\ldots1(2k-2)}(x)}.$$
Now we argue by induction on $k$, using (\ref{reccheb}) and the
fact that $\G_{12}(x)=\frac{1}{1-x}$. Similarly, we get the
explicit formula for $\G_{(k+1)k(k+2)(k-1)(k+3)\ldots1(2k+1)}(x)$.
\end{proof}

Theorem~\ref{th:lislds} and Corollary~\ref{cthf0a} suggest that
there should exist a bijection between the sets
$\MP_n(12\ldots(k+1))$ and
$\MP_n(k(k+1)(k-1)(k+2)(k-2)(k+3)\ldots1(2k))$. Finding it remains
an interesting open question.

\begin{theorem}\label{thfa}
Let $\tau=(\rho',t,k,\theta',1,t-1)\in\MP_k$ such that
$\rho'_a>t>\theta'_b$ for all $a,b$. Let $\rho$ and $\theta$ be
the permutations obtained by decreasing each entry of $\rho'$ by
$t$ and decreasing each entry of $\theta'$ by $1$, respectively.
Then
$$\G_{\tau}(x)=\frac{1-x^2\G_{\rho}(x)\wt\G_{\theta}(x)}{1-x-x^2(\G_{\rho}(x)+\wt\G_{\theta}(x))},$$
where $\wt\G_{\theta}(x)=\frac{1}{1-x-x^2\G_\theta(x)}$.
\end{theorem}
\begin{proof}
By Lemma~\ref{pmp}, we have two possibilities for block
decomposition of a nonempty Motzkin permutation $\pi\in\MP_n$. Let
us write an equation for $\G_{\tau}(x)$. The contribution of the
first decomposition is $x\G_{\tau}(x)$. The second decomposition
contributes $x^2\G_{\rho}(x)\G_{\tau}(x)$ if $\alpha$ avoids
$\rho$, and $x^2(\G_\tau(x)-\G_{\rho}(x))\wt\G_{\theta}(x)$ if
$\alpha$ contains $\rho$. This last case follows from
Theorem~\ref{thf0}, since if $\alpha$ contains $\rho$, $\beta$ has
to avoid $(\theta,1,t-1)$. Hence,
$$\G_{\tau}(x)=1+x\G_\tau(x)+x^2\G_{\rho}(x)\G_{\tau}(x)+x^2(\G_{\tau}(x)-\G_{\rho}(x))\wt\G_{\theta}(x),$$
where $1$ is the contribution of the empty Motzkin permutation.
Solving the above equation we get the desired result.
\end{proof}

For example, for $\tau=546213$ ($\tau=\rho46\theta13$),
Theorem~\ref{thfa} gives
$\G_\tau(x)=\frac{1-2x}{(1-x)(1-2x-x^2)}$.

The last two theorems can be generalized as follows.

\begin{theorem}\label{gthf}
Let
$\tau=(\tau^1,t_1+1,t_0,\tau^2,t_2+1,t_1,\ldots,\tau^m,t_m+1,t_{m-1})$
where $t_{j-1}>\tau^j_a>t_j$ for all $a$ and $j$. We define
$\sigma^j=(\tau^1,t_1+1,t_0,\ldots,\tau^j)$ for $j=2,\ldots,m$,
$\sigma^0=\vr$, and
$\theta^j=(\tau^j,t_j+1,t_{j-1},\ldots,\tau^m,t_m+1,t_{m-1})$ for
$j=1,2,\ldots,m$. Then
$$\G_\tau(x)=1+x\G_\tau(x)+x^2\sum_{j=1}^m(\G_{\sigma^j}(x)-\G_{\sigma^{j-1}})\G_{\theta^j}(x).$$
(By convention, if $\rho$ is a permutation of
\{i+1,i+2,\ldots,i+l\}, then $\G_\rho$ is defined as $\G_{\rho'}$,
where $\rho'$ is obtained from $\rho$ decreasing each entry by
$i$.)
\end{theorem}
\begin{proof}
By Lemma~\ref{pmp}, we have two possibilities for block
decomposition of a nonempty Motzkin permutation $\pi\in\MP_n$. Let
us write an equation for $\G_{\tau}(x)$. The contribution of the
first decomposition is $x\G_{\tau}(x)$. The second decomposition
contributes
$x^2(\G_{\sigma^j}(x)-\G_{\sigma^{j-1}}(x))\G_{\theta^j}(x)$ if
$\alpha$ avoids $\sigma^j$ and contains $\sigma^{j-1}$ (which
happens exactly for one value of $j$), because in this case
$\beta$ must avoid $\theta^j$. Therefore, adding all the
possibilities of contributions with the contribution $1$ for the
empty Motzkin permutation we get the desired result.
\end{proof}

For example, this theorem can be used to obtain the following
result.

\begin{corollary}
{\rm(i)} For all $k\geq3$
$$\G_{(k-1)k12\ldots(k-2)}(x)=\frac{U_{k-3}\ttx}{xU_{k-2}\ttx};$$

{\rm(ii)} For all $k\geq4$
$$\G_{(k-1)(k-2)k12\ldots(k-3)}(x)=\frac{U_{k-4}\ttx-xU_{k-5}\ttx}{x\left(U_{k-3}\ttx-xU_{k-4}\ttx\right)};$$

{\rm(iii)} For all $1\leq t\leq k-3$,
$$\G_{(t+2)(t+3)\ldots(k-1)(t+1)k12\ldots
t}(x)=\frac{U_{k-4}\ttx}{xU_{k-3}\ttx}.$$
\end{corollary}

\subsection{Generalized patterns}\label{ss:gener}
In this section we consider the case of generalized patterns (see
Subsection~\ref{background}), and we study some statistics on
Motzkin permutations.

\subsubsection{Counting occurrences of the generalized patterns
$12\mn3\mn\ldots\mn k$ and $21\mn3\mn\ldots\mn k$} Let
$F(t,X,Y)=F(t,x_2,x_3,\ldots,y_2,y_3,\ldots)$ be the generating
function
$$\sum_{n\geq0}\sum_{\pi\in\MP_n}t^n\prod_{j\geq2}
x_j^{12\mn3\mn\ldots\mn j(\pi)}y_j^{21\mn3\mn\ldots\mn j(\pi)},$$
where $12\mn3\mn\ldots\mn j(\pi)$ and $21\mn3\mn\ldots\mn j(\pi)$
are the number of occurrences of the pattern $12\mn3\mn\ldots\mn
j$ and $21\mn3\mn\ldots\mn j$ in $\pi$, respectively.

\begin{theorem}\label{gp1}
We have
$$F(t,X,Y)=1-\frac{t}{\displaystyle ty_2-\frac{1}
{\displaystyle 1+tx_2(1-y_2y_3)+tx_2y_2y_3F(t,X',Y')}},$$ where
$X'=(x_2x_3,x_3x_4,\ldots)$ and $Y'=(y_2y_3,y_3y_4,\ldots)$. In
other words, the generating function
$F(t,x_2,x_3,\ldots,y_2,y_3,\ldots)$ is given by the continued
fraction
$${\displaystyle 1-\frac{t}
{\displaystyle ty_2-\frac{1}{\displaystyle
1+tx_2-\frac{t^2x_2y_2y_3}{\displaystyle ty_2y_3-\frac{1}
{\displaystyle 1+tx_2x_3-\frac{t^2x_2x_3y_2y_3^2y_4}{\displaystyle
ty_2y_3^2y_4-\frac{1}{\displaystyle
1+tx_2x_3^2x_4-\frac{t^2x_2x_3^2x_4y_2y_3^3y_4^3y_5}{\ddots}}
}}}}}}.$$
\end{theorem}
\begin{proof}
As usual, we consider the two possible block decompositions of a
nonempty Motzkin permutation $\pi\in\MP_n$. Let us write an
equation for $F(t,X,Y)$. The contribution of the first
decomposition is $t+ty_2(F(t,X,Y)-1)$. The contribution of the
second decomposition gives $t^2x_2$, $t^2x_2y_2(F(t,X,Y)-1)$,
$t^2x_2y_2y_3(F(t,X',Y')-1)$, and
$t^2x_2y_2^2y_3(F(t,X,Y)-1)(F(t,X',Y')-1)$ for the four
possibilities (see Figure~\ref{avmp}) $\alpha=\beta=\vr$,
$\alpha=\vr\neq\beta$, $\beta=\vr\neq\alpha$, and
$\beta,\alpha\neq\vr$, respectively. Hence,
$$\begin{array}{l}
F(t,X,Y)=1+t+ty_2(F(t,X,Y)-1)+t^2x_2+t^2x_2y_2y_3(F(t,X'Y')-1)\\
\qquad\qquad\qquad+t^2x_2y_2(F(t,X,Y)-1)
+t^2x_2y_2^2y_3(F(t,X,Y)-1)(F(t,X',Y')-1), \end{array}$$ where $1$
is as usual the contribution of the empty Motzkin permutation.
Simplifying the above equation we get
$$F(t,X,Y)=1-\frac{t}{\displaystyle ty_2-\frac{1}
{\displaystyle 1+tx_2(1-y_2y_3)+tx_2y_2y_3F(t,X',Y')}}.$$ The
second part of the theorem now follows by induction.
\end{proof}

As a corollary of Theorem~\ref{gp1} we recover the distribution of
the number of rises and number of descents on the set of Motzkin
permutations, which also follows easily from Theorem~\ref{th:lds}.

\begin{corollary}
We have
$$\sum_{n\geq0}\sum_{\pi\in\MP_n}t^np^{\#\{\mathrm{rises\ in\ }\pi\}}q^{\#\{\mathrm{descents\ in\ }\pi\}}=
\frac{1-qt-2pq(1-q)t^2-\sqrt{(1-qt)^2-4pqt^2}}{2pq^2t^2}.$$
\end{corollary}

As an application of Theorem~\ref{gp1} let us consider the case of
Motzkin permutations which contain either $12\mn3\mn\ldots\mn k$
or $21\mn3\mn\ldots\mn k$ exactly $r$ times.

\begin{theorem}\label{cgp1a}
Fix $k\geq2$. Let $\G_{12\mn3\mn\ldots\mn k}(x;r)$ be the
generating function for the number of Motzkin permutations which
contain $12\mn3\mn\ldots\mn k$ exactly $r$ times. Then
$$\G_{12\mn3\mn\ldots\mn k}(x;0)=\frac{U_{k-1}\ttx}{xU_k\ttx},$$
and for all $r=1,2,\ldots,k-1$,
$$\G_{12\mn3\mn\ldots\mn k}(x;r)=\frac{x^{r-1}U_{k-2}^{r-1}\ttx}{(1-x)^rU_{k-1}^{r+1}\ttx}.$$
\end{theorem}
\begin{proof}
Let $t=x$, $x_k=y$, $x_j=1$ for all $j\neq k$, and $y_j=1$ for all
$j$. Let $\wt G_k(x,y)$ be the function obtained from $F(t,X,Y)$
after this substitution. Theorem~\ref{gp1} gives
$$\wt G_k(x,y)=1-\dfrac{x}{x-\dfrac{1}{1+x-\dfrac{x^2}{1-\dfrac{x}{1+x-\dfrac{x^2}{\ddots-\dfrac{\ddots}{x-\dfrac{1}{1+xy-\dfrac{x^2y}{x-\dfrac{1}{1+xy^{k+1}-\ddots}}}}}}}}}.$$
Therefore, $\wt G_k(x,y)$ can be expressed as follows. For all
$k\geq2$,
$$\wt G_k(x,y)=1-\frac{x}{x-\dfrac{1}{1+x\wt G_{k-1}(x,y)}},$$ and there exists a
continued fraction $\wt H(x,y)$ such that
$$\wt G_1(x,y)=y-\frac{xy}{x-\dfrac{1}{1+xy^{k+1}\wt H(x,y)}}.$$ Now, using
induction on $k$ together with (\ref{reccheb}) we get that there
exists a formal power series $\wt J(x,y)$ such that
$$\wt G_k(x,y)=\frac{(1-x)U_{k-2}\ttx-xyU_{k-3}\ttx}
{x(1-x)U_{k-1}\ttx-x^2yU_{k-2}\ttx}+y^{k+1}\wt J(x,y).$$ Similarly
as in the proof of Theorem~\ref{thacapp1}, expanding $\wt
G_k(x,y)$ in series about the point $y=0$ gives the desired
result.
\end{proof}

Using the same idea as in Theorem~\ref{cgp1a}, we can apply
Theorem~\ref{gp1} to obtain the following result.

\begin{theorem}\label{cgp1b}
Fix $k\geq2$. Let $\G_{21\mn3\mn\ldots\mn k}(x;r)$ be the
generating function for the number of Motzkin permutations which
contain $21\mn3\mn\ldots\mn k$ exactly $r$ times. Then
$$\G_{21\mn3\mn\ldots\mn k}(x;0)=\frac{U_{k-3}\ttx-xU_{k-4}\ttx}{x\left(U_{k-2}\ttx-xU_{k-3}\ttx\right)},$$
and for all $r=1,2,\ldots,k-1$,
$$\G_{12\mn3\mn\ldots\mn k}(x;r)=\frac{x^{r}(1+x)^rU_{k-2}^{r-1}\ttx}{\left(U_{k-2}\ttx-xU_{k-3}\ttx\right)^{r+1}}.$$
\end{theorem}
%-------------------------------------------------------------------
\small
\subsection*{Acknowledgements}
We would like to thank Marc Noy for helpful comments and
suggestions. The first author was partially supported by a MAE fellowship.

%-------------------------------------------------------------------

\end{document}